\documentclass[12pt]{amsart}
\usepackage{amsthm, amsmath}
\usepackage[psamsfonts]{amssymb} 
\usepackage{amscd}

\newsymbol\rxtimes 226F
\theoremstyle{plain}

\newtheorem{thm}[subsection]{Theorem}
\newtheorem{lem}[subsection]{Lemma}
\newtheorem{prop}[subsection]{Proposition}
\newtheorem{cor}[subsection]{Corollary}

\theoremstyle{remark}

\def\rond{\kern 2pt{\scriptstyle\circ}\kern 2pt}

\def\det{\mathop{\rm det}\nolimits}

\def\Br{\mathop{\rm Br}\nolimits}

\def\Z{\mathbf{Z}}

\def\C{\mathbf{C}}
\def\P{\mathbf{P}}

\def\G{\mathbf{G}}
\def\H{\mathrm{H}}
\def\Mu{\boldsymbol{\mu}}
\def\pg{\mathrm{PGL_2}}
\def\ps{\mathrm{PSL_2}}

\def\so{\mathrm{SO}}
\def\emb{\mathop{\rm Emb}\nolimits _i}
\def\ss{\,{}^\sigma\!}

\def\iso{\vbox{\hbox to .8cm{\hfill{$\scriptstyle\sim$}\hfill}
\nointerlineskip\hbox to .8cm{{\hfill$\longrightarrow $\hfill}} }}

\catcode`\@=11\def\eqalign#1{\null\,\vcenter{\openup\jot\m@th\ialign{
\strut\hfil$\displaystyle{##}$&$\displaystyle{{}##}$
&&\quad\strut$\displaystyle{##}$&$\displaystyle{{}##}$
\crcr#1\crcr}}\,}
\catcode`\@=12

\input amssym.def
\input amssym
\input xy
\xyoption{all}

\begin{document}
\title[Finite subgroups of $\pg(K)$]{Finite subgroups of $\mathbf{PGL_2(\boldsymbol{K})}$}
\author[Arnaud Beauville]{Arnaud Beauville}
\address{Laboratoire J.-A. Dieudonn\'e\\
UMR 6621 du CNRS\\
Universit\'e de Nice\\
Parc Valrose\\
F-06108 Nice cedex 2, France}
\email{arnaud.beauville@unice.fr}
% \dedicatory{To Ramanan on  his 70th birthday}
\date{\today}
 
\begin{abstract}
We classify, up to conjugacy, the finite subgroups of $\pg(K)$ of order prime to $\mathrm{char}(K)$.
\end{abstract}

\maketitle

\section*{Introduction}
\par The aim of this note is to describe, up to conjugacy, the finite subgroups of    $\pg(K)$, for an arbitrary field $K$. Throughout the paper, \emph{we consider only subgroups whose order is prime to the characteristic of} $K$. 
\par When $K=\C$, or more generally when $K$ is algebraically closed, the answer is well known: any such group is isomorphic to $\Z/r$,  $D_r$ (the dihedral group), 
$\mathfrak{A}_4$, $\mathfrak{S}_4$ or $\mathfrak{A}_5$, and  there is only one conjugacy class for each of these groups. If $K$ is arbitrary, the group $\pg(K)$ is contained in $\pg(\overline{K})$, so the subgroups of $\pg(K)$ are among the previous list; it is not difficult to decide which subgroups occur for a given field $K$, see \S 1.
\par So the only question left is to describe the conjugacy classes in $\pg(K)$ of the subgroups in the list.  In \S 2 we give a general answer for subgroups of $G(K)$, for an algebraic group $G$, in terms of (non-abelian) Galois cohomology. We illustrate the method on one example in \S 3, and  apply it to the case $G=\pg$ in \S 4.
\par The motivation for looking at this question was to understand the appearance of the Brauer group in  the case of $(\Z/2)^2$ considered in \cite{B}. The result is somewhat disappointing, as it turns out that this case (which could be treated directly, as in \cite{B})  is the only one where some second Galois cohomology group plays a role.  At least our method explains this role, and hopefully may be useful in other situations.
\medskip
\section{The possible subgroups}
\par   We repeat that whenever we mention a finite group, we always assume that its order is prime to the characteristic of $K$. 
The following is classical (see \cite{S2}, 2.5).
\begin{prop} \label{list}
$1)$ $\pg(K)$ contains $\Z/r$ and $D_r\,$\footnote{We denote by $D_r$ the dihedral group with $2r$ elements.} if and only if $K$ contains $\zeta +\zeta ^{-1}$ for some  primitive $r$-th root of unity $\zeta $.
\par  $ 2)$ $\pg(K)$ contains $\mathfrak{ A}_4$ and $\mathfrak{ S}_4$ if and only $-1$ is the sum of two squares in $K$.
%if there exists $\alpha ,\beta $ in $K$ with $\alpha ^2+\beta ^2=-1$.
\par   $3)$ $\pg(K)$ contains $\mathfrak{ A}_5$  if and only if $-1$ is the sum of two squares  and $5$ is a square in $K$.
\end{prop}

\noindent\textit{Proof} : One way to prove this is to use the isomorphism $\pg(K)\iso$ $\so(K,q)$, where $q$ is the  quadratic form $q(x,y,z)=x^2+yz$ on $K^3$ (\cite{D}, II.9). If a group $H$ embeds into $\so(K,q)$, we have a faithful representation $\rho $ of $H$ in $K^3$, which preserves an indefinite quadratic form. 

\smallskip
\par   $\bullet\ $ Case $H=\Z/r$ : let $g$ be a generator; the existence of $q$ forces  the eigenvalues of $\rho (g)$ in $\overline{K}$ to be of the form $(\zeta,  \zeta ^{-1}, 1)$, with $\zeta $ a primitive $r$-th root of 1. This implies $\zeta +\zeta ^{-1}\in K$. Conversely, if $\lambda  :=\zeta +\zeta ^{-1}$ is in $K$, the homography $\displaystyle z\mapsto {(\lambda +1)z-1\over  z+1}$ is an element of order $r$ of $\pg(K)$. 
\smallskip
\par   $\bullet\ $ Case $H=D_r$ : by the previous case, if $D_r\subset \pg(K)$, $\lambda  :=\zeta +\zeta ^{-1}$ is in $K$. Conversely if $\lambda  \in K$, the homographies ${z\mapsto 1/z}$ and  $ \displaystyle z\mapsto {(\lambda +1)z-1\over  z+1}$  generate a subgroup of 
$\pg(K)$ isomorphic to $D_r$.
\medskip
\par   $\bullet\ $ Cases $H=\mathfrak{ A}_4, \mathfrak{ S}_4$ or $\mathfrak{ A}_5\,$. 
The representation $\rho $ must be irreducible. 
Each of the  groups $\mathfrak{ A}_4$ and $\mathfrak{ S}_4$ has exactly one  irreducible 3-dimensional representation with trivial determinant, which 
is defined over the prime field; the only invariant quadratic form (up to a scalar) is the standard form $q_0(x,y,z)=x^2+y^2+z^2$. Thus 
$\mathfrak{A}_4$ and $\mathfrak{S}_4$  are contained in $\pg(K)$ if and only if $q_0$ is equivalent to $\lambda q$ for some $\lambda \in K^*$, which means that $q_0$ represents 0.
\par   Since $\mathfrak{ A}_5$ contains  elements of order 5,  the condition $\sqrt{5}\in K$ is necessary. Suppose this is the case, and put $\varphi =\frac{1}{2}(1+\sqrt{5} )$; the subgroup of $\so(K,q_0)$ preserving the
icosahedron with vertices  $$\{(\pm 1,0,\pm \varphi )\ ,\ (\pm \varphi ,\pm 1,0 )\ ,\ (0,\pm \varphi,\pm 1 )\}$$ is isomorphic to $\mathfrak{ A}_5$. It follows as above that $\mathfrak{ A}_5$ embeds in $\so(K,q)$ if and only if $q_0$ represents 0.\qed

\medskip

\section{Some Galois cohomology}
\subsection{}  In this section we consider an algebraic group $G$ over $K$, and a subgroup $H\subset G(K)$. We choose a separable closure $K_s$ of $K$, and put $\mathfrak{g}:=\mathrm{Gal}(K_s/K)$. We are interested in the set of embeddings $H \hookrightarrow G(K)$ which are conjugate in $G(K_s)$ to the natural inclusion $i:H \hookrightarrow G(K)$, modulo conjugacy by an element of $G(K)$.  We denote this (pointed) set by $\emb(H,G(K))$. 
\par  We will use the standard conventions for non-abelian cohomology, as explained for instance in \cite{S3}, ch. I, \S 5. We will also use the notation of \cite{S3} for Galois cohomology: if $G$ is an algebraic group over $K$, we put $\H^i(K,G):=\H^i(\mathfrak{g},G(K_s))$.

\begin{prop}\label{co}
Let $Z$ be the centralizer of $H$ in $G(K_s)$. The pointed set $\emb(H,G(K))$ is canonically isomorphic to the kernel of the natu\-ral map
$\H^1(K,Z)\rightarrow \H^1(K,G)$.
\end{prop}
\noindent\textit{Proof} : Let $X\subset G(K_s)$ be the subset   of elements $g$ such that \allowbreak $g^{-1}\ss g\in Z$ for all $\sigma \in \mathfrak{g}$. The group $G(K)$ (resp. $Z$) acts on $X$ by left (resp. right) multiplication. By \cite{S3}, ch. I, 5.4, cor. 1, the kernel of $\H^1(K,Z)\rightarrow \H^1(K,G)$ is identified with the (left) quotient by $G(K)$ of  the subset of $\mathfrak{g}$-invariant elements in $G(K_s)/Z$; but this subset is by definition  $X/Z$, so we can identify our kernel to the double quotient $G(K)\backslash X/Z$. 
\par   For every $g\in X$, the conjugate embedding $gig^{-1}$  belongs to \break $\emb(H,G(K))$. Any element  
 $j\in \emb(H,G(K))$ is of the form $gig^{-1}$ for some $g\in G(K_s)$; for  $\sigma \in \mathfrak{g} $, the element $\ss g$ again conjugates $i$ to $j$, hence $g^{-1}\ss g\in Z$ and $g\in X$. Thus the map $g\mapsto gig^{-1}$ from $X$ to $\emb(H,G(K))$ is surjective. Two elements $g$ and $g'$ of $X$ give the same element in $\emb(H,G(K))$ if and only if $g'$ belongs to the double coset $G(K)gZ$. Therefore the above map induces a canonical bijection $G(K)\backslash X/Z\iso\emb(H,G(K))$.\qed
 \medskip
 \def\co{\mathop{\rm Conj}\nolimits (H,G(K))}
\subsection{} Let us write down the correspondence explicitly:  a class in our kernel is represented by a 1-cocycle $\mathfrak{g}\rightarrow Z$ which becomes a coboundary in $G$, hence is of the form $\sigma \mapsto g^{-1}\ss g$ for some $g\in X$; we associate to this class the embedding $gig^{-1}$.\label{corr}
  \subsection{}   We are actually more interested in  the set   $\co$  of subgroups of $G(K)$ which are conjugate to $H$ in $G(K_s)$, modulo conjugacy by $G(K)$. Associating to an embedding its image defines a surjective map $im:\emb(H,G(K))\rightarrow \co$. 
The normalizer $N$ of $H$ in $G(K_s)$ acts on $H$ by automor\-phisms, hence also on $\emb(H,G(K))$. Two embeddings with the same image differ by an automorphism of $H$, 
which must be induced by an element of $N$ if the embeddings are conjugate under $G(K_s)$. It follows that $im$ \emph{induces an isomorphism} $\emb(H,G(K))/N\iso \co$.

\subsection{} Let us translate this in cohomological terms. Let $\H^1(K,Z)_0$ denote the kernel of the map $\H^1(K,Z)\rightarrow \H^1(K,G)$.
  An element $n$ of $N$ acts on $\emb(H,G(K))$ by $j\mapsto j\rond \mathrm{int}(n^{-1})$; if $j=gig^{-1}$, this amounts to replace $g$ by $gn$, hence the 1-cocycle $\varphi :\sigma \mapsto g^{-1}\ss g$ by $n^{-1}\varphi \ss n$. This formula defines an action of $N$ on $\H^1(K,Z)$ which preserves  $\H^1(K,Z)_0$;  \emph{the map  $g\mapsto gHg^{-1}$ induces an isomorphism of pointed sets} 
$\H^1(K,Z)_0/N \iso \co$.
\label{N}
\medskip
\section{An example}
\subsection{}  In this section we fix an integer $r\geq 2$, prime to $\mathrm{char}(K)$, and we  assume that $K$ contains a primitive $r$-th root of unity $\zeta $. We consider the matrices $A,B\in\mathrm{M}_r(K)$ defined on the canonical basis $(e_1,\ldots ,e_r)$ of $K^r$ by
$$A\cdot e_i=e_{i+1}\quad ,\quad B\cdot e_i=\zeta ^ie_i$$
for $1\leq i\leq r$, with the convention $e_{r+1}=e_1$.
\par  The matrices $A$ and $B$ generate the $K$-algebra $\mathrm{M}_r(K)$, with the relations
$$A^r=B^r=I\quad ,\quad BA=\zeta AB\ .$$ 
Their classes $\bar A,\bar B$ in $\mathrm{PGL}_r(K)$ commute; we consider the embedding $i:(\Z/r)^2\hookrightarrow \mathrm{PGL}_r(K)$ which maps the two basis vectors to $\bar A$ and $\bar B$.
The image $H$ of $i$ is its own centralizer; in particular, $H$ is a maximal commutative subgroup of  $\mathrm{PGL}_r(K)$. %When $r$ is prime it is the only $r$-subgroup of $\mathrm{PGL}_r(K)$ not contained in a maximal torus (\cite{Bo}, 6.4).
\par  By the Kummer exact sequence (and the choice of $\zeta $), the group $\H^1(K,\Z/r)$ is identified with $K^*/K^{*r}$; the pointed set $\H^1(K,\mathrm{PGL}_r)$ can be viewed as the set of isomorphism classes of central simple $K$-algebras of dimension $r^2$ (\cite{S1}, X.5).
\begin{lem}
Let $\alpha ,\beta \in K^*$, and let $\bar \alpha ,\bar \beta $ be their images in $K^*/K^{*r}$. The map $\H^1(i):\H^1(K,\Z/r)^2\rightarrow \H^1(K,\mathrm{PGL}_r)$ associates to $(\bar \alpha ,\bar\beta )$ the class of the cyclic  $K$-algebra $A_{\alpha ,\beta}$ generated  by two variables $x,y$ with the relations $x^r=\alpha $, $y^r=\beta $, $yx=\zeta xy$.
\end{lem}
\noindent\textit{Proof} : We choose $\alpha ',\beta '$ in $K_s$ with $\alpha '^r=\alpha $ and $\beta '^r=\beta $. The Kummer isomorphism associates to $(\alpha ,\beta )$ the homomorphism  $(a,b) :\mathfrak{g}\rightarrow (\Z/r)^2 $ defined by
$$\ss\alpha '=\zeta ^{a(\sigma )}\alpha '\qquad \ss\beta '=\zeta ^{b(\sigma )}\beta '\qquad\hbox{for each }\sigma \in\mathfrak{g}\ .$$ 
Its image in $\H^1(K,\mathrm{PGL}_r(K_s))$ is the class of the 1-cocycle $\sigma \mapsto \bar A^{a(\sigma )} \bar B^{b(\sigma )}$.
\par  Now let us recall how we
associate to the algebra $A_{\alpha ,\beta}$ a cohomology class  $[A_{\alpha ,\beta }]$ in $ \H^1(K,\mathrm{PGL}_r)$ (\textit{loc. cit.}). We choose an isomorphism of $K_s$-algebras $u:\mathrm{M}_r(K_s)\iso A_{\alpha ,\beta }\otimes_KK_s $.
For each $\sigma \in\mathfrak{g}$, $u^{-1}\ss u$ is an automorphism of $\mathrm{M}_r(K_s)$, hence of the form $\mathrm{int}(g_\sigma )$ for some $g_\sigma $ in $\mathrm{PGL}_r(K_s)$. 
Then $[A_{\alpha ,\beta }]$ is the class of the 1-cocycle $\sigma \mapsto g_\sigma $. 

\par  In our case we define $u$ on the generators $A,B$ by $u(A)=\beta 'y^{-1}$, $u(B)=\alpha  '^{-1}x$. Then the automorphism 
$u^{-1}\ss u$ multiplies $A$ by $\zeta ^{b(\sigma )}$ and $B$ by $\zeta ^{-a(\sigma )}$, which gives $g_\sigma=\bar A^{a(\sigma )} \bar B^{b(\sigma )}$ as above.\qed

\subsection{}\label{symbol}  The exact sequence
$$1\rightarrow \G_m\rightarrow \mathrm{GL}_r\rightarrow \mathrm{PGL}_r\rightarrow 1$$gives rise to a coboundary homomorphism $\partial_r:\H^1(K,\mathrm{PGL}_r)\rightarrow$\break $\H^2(K,\G_m)=\Br(K)$ which is injective (\textit{loc. cit.}). The class $\partial_r[A_{\alpha ,\beta }]\in\Br(K)$ is the \emph{symbol} $(\alpha ,\beta )_r$; it depends only on the classes of $\alpha $ and $\beta $ (mod. $K^{*r}$). 
The map $(\ ,\ )_r:(K^*/K^{*r})^2\rightarrow \Br(K)$ is bilinear and alternating. Since $\partial_r$ is injective, we find:
\begin{prop}\label{ex}
The set $\emb((\Z/r)^2, \mathrm{PGL}_r(K))$ is isomorphic to the set of couples $(\alpha ,\beta )$ in $(K^*/K^{*r})^2$ such that $(\alpha ,\beta )_r=0$.\qed
\end{prop}
\par  We will describe the correspondence more explicitely in the case $r=2$ in the next section.
\medskip
\section{Conjugacy classes in $\pg(K)$}
  
\begin{prop}\label{unique}
 Assume that $K$ is separably closed. Two finite subgroups of $\pg(K)$ which are isomorphic $($and of order prime to \allowbreak$\mathrm{char}(K))$ are conjugate. 
\end{prop}
\noindent\textit{Proof }: Again this is certainly well-known; we give a quick proof for completeness. The possible subgroups are those which appear in  Proposition \ref{list}. 
\par An element of order $r$ of $\pg(K)$ comes from a diagonalizable element of $\mathrm{GL}_2(K)$, hence is conjugate to the homothety $z\mapsto \zeta z$ for some $\zeta \in \Mu _r(K)\,$\footnote{As usual we denote by $\Mu _r(K)$ the group of $r$-th roots of unity in $K$.}; thus a cyclic subgroup of order $r$ of $\pg(K)$ is conjugate to the group $H_r$
of homotheties $z\mapsto \lambda z,\ \lambda \in \Mu _r(K)$. 
\par There is only one group $D_r$ containing $H_r$, namely the subgroup  generated by $H_r$ and the involution $z\mapsto 1/z$; it follows that all dihedral subgroups  of order $2r$ are conjugate to this subgroup.
\par For the three remaining groups, we use again the isomorphism \allowbreak $\pg(K)\iso \so_3(K)$. The groups $\mathfrak{A}_4$ and $\mathfrak{S}_4$ have exactly one irreducible representation of dimension 3 with trivial determinant, while $\mathfrak{A}_5$  has two such representations which differ by an outer automorphism: this is elementary in characteristic 0, and the general case follows by \cite{I}, ch. 15. 
Therefore two isomorphic subgroups $H$ and $H'$ of $\so_3(K)$ of this type are conjugate in $\mathrm{GL}_3(K)$.
The only quadratic forms preserved by $H$ or $H'$ are the multiple of the standard form; thus the element $g$ of  $\mathrm{GL}_3(K)$ which conjugates $H$ to $H'$ must satisfy ${}^t\!g\,g=\lambda I$ for some $\lambda \in K$. Replacing $g$ by $\pm\mu g$, with $\mu ^2=\lambda ^{-1}$, we have $g\in\so_3(K)$, hence our assertion.\qed
\medskip
\par Recall that the determinant induces a homomorphism  $\ \overline{\det}:$ \break $\pg(K)\rightarrow K^*/K^{*2}$.

\begin{thm}
 $1)\ \pg(K)$ contains only one conjugacy class of subgroups isomorphic to $\Z/r\ (r>2)$,  $\mathfrak{ A}_4$,  $\mathfrak{ S}_4$ or $\mathfrak{ A}_5$.
\par $2)$ The conjugacy classes of cyclic subgroups of order $2$ of $\pg(K)$ are parametrized by $K^*/K^{*2}$: to $\alpha \in K^*\ \mathrm{(mod.}\, K^{*2})$ corresponds the involution $z \mapsto \alpha/ z$.
\par $3)$ The homomorphism $\overline{\det}:\pg(K)\rightarrow K^*/K^{*2}$ induces a bijective correspondence between:
\par  $\bullet$ conjugacy classes of subgroups of $\pg(K)$ isomorphic to $(\Z/2)^2$;
\par  $\bullet$ subgroups $G\subset K^*/K^{*2}$ of order $\leq 4$, such that $(-\alpha ,-\beta )_2=0$ for all $\alpha ,\beta $ in $G$ \hbox{\rm (see (\ref{symbol})).}

\par $4)$ Assume that $\Mu _r(K)$ has order $r$. The conjugacy classes of  subgroups $D_r$ of $\pg(K)$ are  parametrized by  $K^*/K^{*2}\Mu _r(K)$. The subgroup corresponding to $\alpha \in K^*\  {\rm (mod.}\ K^{*2}\Mu _r(K)\,)$ consists of the homo\-graphies $z\mapsto \zeta  z$ and $z\mapsto \alpha\eta  / z$, for $\zeta ,\eta  \in \Mu _r(K)$.
\end{thm}

\noindent\textit{Proof} : Using Proposition \ref{unique} we can apply the method of \S 3. 
 We give the list of the subgroups of $\pg(K_s)$ and their centralizers:
 \medskip
\def\tvi{\vrule height 16pt depth 8pt width 0pt}
\def\tv{\tvi\vrule}
\def\tvj{\tvi\vrule width 1pt}
\def\n{\noalign{\hrule}}
\def\h{\hfill\kern5pt }
$$\hss\vbox{\offinterlineskip
\halign{\tv$#$&\h$#$\h\tvj &\h\ $#$\h\ \tv &\h $#$\h \tv& &\h $#$\h \tv &\h $#$\h \tv&\h $#$\h \tv &\h $#$\h \tv\cr\n
&H&\Z/2&\Z/r \ (r>2) &\Z/2\times \Z/2& D_r \ (r>2)& \mathfrak{ A}_4 & \mathfrak{ S}_4 & \mathfrak{ A}_5\cr\n
&Z&\G_m \rxtimes  \Z/2 & \G_m & \Z/2\times \Z/2 & \Z/2 & 1& 1&1 \cr\n
}}\hss$$

\par  In case 1), we have $\H^1(K,Z)=\{1\}$  (using $\H^1(K,\G_m)$ $=\{1\}$). The result follows from  (\ref{N}).   \smallskip
\par Case 2): We apply Proposition \ref{co}, taking for $i(\bar 1)$ the involution $z\mapsto 1/z$. The centralizer $Z$ is the semi-direct product of a torus $\G_m $ and the subgroup $\Z/2$  generated by the involution $\iota :z\mapsto -z$. The pointed set $\H^1(K,\G_m \rxtimes \Z/2)$ is identified with $\H^1(K, \Z/2)\cong K^*/K^{*2}$. The map
$\H^ 1(K, \Z/2)\rightarrow \H^ 1(K, \pg)$ is trivial, for instance because the injection $\Z/2\rightarrow \pg$ factors through a torus $\G_m$. Hence the set of conjugacy classes of involutions of $\pg(K)$ is identified with $K^*/K^{*2}$.
\par To describe the correspondence explicitely we follow (\ref{corr}). Let $\alpha \in K^*$, and let $\alpha ' \in K_s^*$ such that $\alpha'^2=\alpha $; the class of $\alpha $ (mod. $K^{*2}$) corresponds to the class of the 1-cocycle
 $a :\mathfrak{g}\rightarrow \Z/2$ given by $ \ss\alpha ' =(-1)^{a (\sigma )}\alpha'$. In $\pg(K_s)$ we have $i(a (\sigma ))=g^{-1}\ss g$, where $g$ is the homography $z\mapsto \alpha ' z$. Thus the subgroup corresponding to $\alpha $ is $gH g^{-1} $, which is generated by the involution $z\mapsto \alpha /z$.
\medskip
\par Case 3): Let $i:(\Z/2)^2\hookrightarrow \pg(K)$ be the embedding which maps the  basis vectors $e_1$ and $e_2$ to the involutions $z\mapsto 1/z$ and $z\mapsto -z$. By Proposition \ref{ex} the set ${\mathop{\rm Emb}\nolimits _i((\Z/2)^2
,\pg(K))}$ is canonically identified to the set of couples $(\alpha ,\beta )$ in $(K^*/K^{*2})^2$ with $(\alpha ,\beta )_2=0$.

\par Again we make the correspondence explicit following  (\ref{corr}). Let $\alpha ,\beta \in K^*$ with $(\alpha ,\beta )_2=0$. This means that the conic $x^2-\alpha y^2-\beta z^2=0$ is isomorphic to $\P^1_K$, thus there exists $\lambda ,\mu $ in $K$ with $\lambda^2 -\alpha -\beta \mu ^2=0$. We choose $\alpha '$ and $\beta '$ in $K_s$ such that $\alpha '^2=\alpha $ and $\beta '^2=\beta $; as above we define the homomorphisms $a $ and $b :\mathfrak{g}\rightarrow \Z/2$  by 
$$\ss\alpha ' =(-1)^{a (\sigma )}\alpha' \quad\hbox{and}\quad \ss\beta ' =(-1)^{b (\sigma )}\beta '\quad\hbox{for each }\sigma \in\mathfrak{g}\ .$$
Put $\displaystyle\theta := {\beta '\mu \over \lambda +\alpha '}={\lambda -\alpha '\over \beta '\mu }$; let $g\in \pg(K_s)$ be the homography $\displaystyle z\mapsto \alpha '{z-\theta \over z+\theta }$ .  An easy computation gives
$$g^{-1}\ss g = i(a (\sigma ),b(\sigma ))\ .$$
Thus the embedding of $(\Z/2)^2$ associated to $(\alpha ,\beta )$ is $gig^{-1}$; it maps $e_1$ to the homography $\displaystyle h_1:z\mapsto {\lambda u-\alpha \over z-\lambda }$ , and $e_2$ to $h_2:z\mapsto \alpha /z$ . Note that $\overline{\det}(h_1)=-\beta $ and $\overline{\det}(h_2)=-\alpha $. 
\par Now we have to take into account the action of the normalizer $N$ of $H$ in $\pg(K_s)$. This is the subgroup $\mathfrak{ S}_4$ generated by $H$ and the homographies
$$n_1:z\mapsto {z+1\over z-1}\quad , \quad n_2:z\mapsto \iota z\ ,$$where $\iota $ is a square root of $-1$.
We apply the recipe of (\ref{N}).   Since $n_1\in \pg(K)$, it acts on $\H^1(K,H)$ through its action on $H$, which permutes $e_1$ and $e_2$; thus it maps $(\alpha ,\beta )\in (K^*/K^{*2})\times (K^*/K^{*2})$ to $(\beta ,\alpha )$. The action of $n_2$  on $H$  fixes $e_2$ and exchanges $e_1$ with $e_1+e_2$; to get the action on $\H^1(K,H)$ we have to multiply by the class of the cocycle $\sigma \mapsto n_2^{-1}\ss n_2$, that is, $\sigma \mapsto i\bigl((\sigma (\iota )/\iota )\,e_2\bigr)$. Hence $n_2$ acts on $\H^1(K,H)$ by
$$n_2\cdot (\alpha ,\beta )=(\alpha  ,-\alpha \beta )\ .$$
\par  Let $G_{\alpha ,\beta }$ be the subgroup of $K^*/K^{*2}$ generated by $-\alpha $ and $-\beta $; it is the image of $H$ by the homomorphism $\overline{\det}:\pg(K)\rightarrow K^*/K^{*2}$.  If $G_{\alpha ,\beta }\cong (\Z/2)^2$, the orbit $N\cdot (\alpha ,\beta )$ in $(K^*/K^{*2})\times (K^*/K^{*2})$ has 6 elements, which are the couples $(-x,-y)$ with $x,y\in G_{\alpha ,\beta }$, $x\neq y$.  If $G_{\alpha ,\beta }\cong (\Z/2)$, the orbit has 3 elements, which are the couples $(-x,-y)$ with $x,y\in G_{\alpha ,\beta }$, $(x,y)\neq (1,1)$. Finally if  $G_{\alpha ,\beta }$ is trivial the orbit consists only of $(-1,-1)$. Thus the  conjugacy classes of subgroups $(\Z/2)^2$ in $\pg(K)$ are parametrized   by the subgroups  $G\subset K^*/K^{*2}$ of order $\leq 4$, with the property $(-\alpha ,-\beta )_2=0$ for each $\alpha ,\beta $ in $G$.
\medskip
\par Case 4):  The group $D_r$ is generated by two elements $s,t$ with the relations $s^2=t^r=1$ and $sts=t^{-1}$. We choose a primitive $r$-th root of unity $\zeta $ and consider the embedding $i:D_r\hookrightarrow \pg(K)$ such that $i(s)$ is 
 the involution $z\mapsto 1/z$ and $i(t)$ the homothety $z\mapsto \zeta z$. The centralizer is $\Z/2$, generated by the involution $z\mapsto -z$. As in case 2) it follows 
that  $\emb(D_r, \pg(K))$ is isomorphic to $\H^ 1(K,\Z/2)$. Also the previous argument shows that the embedding corresponding to $\alpha \in K^*$ is the conjugate of $i$ by the homography $z\mapsto\alpha ' z$, with $\alpha '^2=\alpha $, so it maps $s$ to $z\mapsto \alpha /z$ and $t$ to $z\mapsto \zeta z$. 
\par To complete the picture we have to take into account the action of the normalizer $N$ of $i(D_r)$ in $\pg(K_s)$. This is the subgroup $D_{2r}$ generated by $i(s):z\mapsto 1/z$ and the homothety $n:z\mapsto \eta z$, where $\eta \in K_s$ is a primitive $2r$-th root of unity. The action 
of $i(s)$ is trivial, and $n$ acts by multiplication by the cocycle $\sigma \mapsto n ^{-1}\ss n $, which corresponds to the class of $\eta ^2 $ in $K^*/K^{*2}$. Since $\eta ^2$ generates $\Mu _r(K)$, the assertion 4) follows.\qed

\bigskip
%\section*{Acknowledgements}

\end{document}